\documentclass{article}
\pdfoutput=1
\usepackage{amsmath}
\usepackage{amsfonts}
\usepackage{amssymb}
\usepackage{amsthm}
\usepackage[pdftex]{graphicx}
\usepackage{epstopdf}
\usepackage{cite}
\usepackage{makecell}
\usepackage{bm}%
\providecommand{\U}[1]{\protect\rule{.1in}{.1in}}
\providecommand{\U}[1]{\protect\rule{.1in}{.1in}}
\usepackage[margin=0.8in]{geometry}
\newtheorem{theorem}{Theorem}

\newtheorem{definition}[theorem]{Definition}

\newtheorem{proposition}[theorem]{Proposition}

\begin{document}

	\title{\LARGE \textbf{Exact recursive updating of uncertainty sets for discrete-time plants with a lag}}
	\author{Robin Hill$^{1}$, Yousong Luo$^{2}$  \thanks{$^{1}%
			$Department of Electrical and Electronic Engineering, University of Melbourne,
			Melbourne, Vic 3010, Australia \texttt{{\small robin.hill@unimelb.edu.au}}%
		}\thanks{$^{2}$School of Science, RMIT
		University, 124 Latrobe St, Melbourne, 3001, Australia
		\texttt{{\small yluo@rmit.edu.au}}}}
\maketitle

\begin{abstract}
In \cite{Hill-luo-Schw-2016} there are new results concerning the polytopic set of possible
states of a linear discrete-time SISO system subject to bounded disturbances from
measurements corrupted by bounded noise. Using these results we construct an algorithm which, for the special case of a plant with a lag, recursively updates these polytopic sets when new measurements arrive. 
 
\end{abstract}

\pagestyle{plain} 

\section{INTRODUCTION}
\setcounter{page}{1}

Determining the possible states of a dynamic system from noisy measurements, when the input and measurement errors are unknown
except for bounds on their magnitude, is a computationally challenging task. We use the term \textit{uncertainty set} to describe such a  set of states consistent with past measurements and \textit{a priori} assumptions on the exact plant model and bounds on the input and measurement noise. 

The problem of recursively calculating uncertainty sets was formulated, and in principle solved, in the seminal papers \cite{Witsenhausen-poss-states-1968,Schweppe-1968,Bertsekas-etal-1971}. There is an extensive literature on
the topic of uncertainty set propagation. See \cite{Fogel-Huang-82} and \cite{Ninness_Goodwin-95} for background, the survey paper \cite{Milanese-Vicino-91} and the book \cite{blanchini_miani}. Some of the many other papers which consider this problem are \cite{Blanchini-Sznaier-2012}, \cite{Rosa_silvestre_athans-2014}, \cite{Stoorvogel-1996} and \cite{Tempo-1988}.  In real-time applications, for example fault detection and isolation \cite{Rosa-CDC2010}, existing exact algorithms are of limited use because of their computational complexity. For this reason there has been a lot of research recently on the use of zonotopes and constrained zonotopes to approximate the exact polytopic uncertainty set, for example \cite{Combastel2015265,Alamo-etal-2008,Scott2016126}. 

In \cite{Hill-luo-Schw-2016} two theorems have been presented from which new algorithms can be derived which update polytopic uncertainty sets exactly for linear, SISO, discrete-time plants. We believe that these algorithms run faster than any current exact schemes. The particular form that the new algorithms take depends on properties of the plant, and in this paper we present details for the special case of plants without a feedthrough term, that is plants with a lag. The theorems quoted from \cite{Hill-luo-Schw-2016} cover both the lag and non-lag cases, but the details of algorithm development for a plant with a lag are a lot simpler than the non-lag case treated in \cite{Hill-luo-Schw-2016}.

Exact, recursive methods currently require polytopic projection techniques, for example Fourier-Motzkin elimination, see \cite{Keerthi_Gilbert-1987}, \cite{Rakoviv-Mayne-2004} and \cite{Shamma_Kuang-1999}. In these papers it is the identification of redundant inequality constraints that is most demanding computationally. The algorithm we present relies on updating the vertex-facet incidence matrix and does not generate any redundant constraints. Another advantage is that, along with the incidence matrix, both the vertices and facets of the uncertainty set are updated. Existing methods update either the facets or vertices, but not both. It is well known that converting from a representation of the polytope in terms of inequality constraints to a vertex representation is not easy, so having both at our disposal is valuable.

We give representative plots which demonstrate that use our algorithm requires a smaller update time than two other commonly  used  methods.
\section{Basic Setup}
The material in this section and the next is a summary of results in \cite{Hill-luo-Schw-2016}, specialized to the case of a plant with no direct feedthrough. The plant $P$, a linear, time-invariant, causal discrete-time, $m^{\rm th}$
order scalar system, is assumed known. There are two sources of
uncertainty, an input noise disturbance $(u_{j})_{j=0}^{\infty}=\mathbf{u},$
and output measurement noise $(w_{j})_{j=0}^{\infty}=\mathbf{w}.$ The plant
output is $(y_{j})_{j=0}^{\infty}=\mathbf{y}$, and the measurement at time $k$ 
is $z_{k}=y_{k}+w_{k}.$ The initial state, at time $k=0,$ is assumed to be
known exactly, but nothing is known about the uncertainties except that they
satisfy $\left\vert u_{j}\right\vert \leq1$ and $\left\vert w_{j}\right\vert
\leq1.$ 

Given an initial state $\mathbf{x}_{0},$ the measurement history $z_{0},z_{1}%
,\ldots,z_{k},$ and the plant dynamics, we seek the uncertainty set
at time $k+1,$ denoted $S_{k+1};$ it is the set of possible states consistent with
the measurements up to and including $z_{k}$, and is easily seen to be a closed, convex polytope.

\subsection{Notation\label{sectnotprelim}}
The names of incidence vectors and matrices, whose elements are zeros or ones, will start with the letter $\mathcal{I}$. The $\lambda
$-transform (generating function) of an arbitrary sequence $\mathbf{y}%
=(y_{k})_{k=0}^{\infty}$ is defined to be $\hat{\mathbf{y}}(\lambda):=\sum
_{k=0}^{\infty}y_{k}\lambda^{k}.$  Real Euclidean space of dimension $m$ is denoted $\mathbb{R}^{m}$, where $m$ is the order of the plant $P$. States of the plant $P$ are represented by vectors, or points, in  $\mathbb{R}^{m}$. 

\subsection{Transfer function description}

The plant for the estimation system has the transfer function representation
$P(\lambda)=\hat{\mathbf{n}}(\lambda)/\hat{\mathbf{d}}(\lambda)$ where
\begin{align}
\hat{\mathbf{n}}(\lambda)  &  =n_{0}+n_{1}\lambda+n_{2}\lambda^{2}%
+\dotsb+n_{m}\lambda^{m}\nonumber\\ 
\hat{\mathbf{d}}(\lambda)  &  =d_0+d_{1}\lambda+d_{2}\lambda^{2}+\dotsb
+d_{m}\lambda^{m},\nonumber
\end{align}
$m\geq1$, $\hat{\mathbf{n}}(\lambda)$ and $\hat{\mathbf{d}%
}(\lambda)$ are coprime polynomials with real coefficients, and
it is assumed that both the plant $P(\lambda)$ and the plant $P^{\ast}%
(\lambda)$ for the dual system, defined below, are causal, implying
$d_{0}\neq0$ and $d_{m}\neq0$. It is assumed that the plant has a lag so $n_0=0$, and without loss of generality we take $d_{0}=1.$ 

Assuming zero initial conditions, $\mathbf{y}$ and $\mathbf{u}$ are related
by $\hat{\mathbf{d}}(\lambda)\hat{\mathbf{y}}(\lambda)=\hat{\mathbf{n}}(\lambda)\hat{\mathbf{u}}(\lambda).$ 
The indeterminate $\lambda$ represents the unit delay operator in the time domain.
\subsection{State-space representations\label{sssection}}

A state-space description of the estimation system is
\begin{align}
\mathbf{x}_{k+1}  &  =\mathbf{A}\mathbf{x}_{k}+\mathbf{B}u_{k}\label{primalstate1}\\
y_{k}  &  =\mathbf{C}\mathbf{x}_{k}\label{primalstate2}\\
z_{k}  &  =y_{k}+w_{k}\nonumber
\end{align}

where
\begin{subequations}
	\begin{align}
	\mathbf{A}  &  =\left[
	\begin{array}
	[c]{cc}%
	\mathbf{0} & \mathbf{I}_{m-1}\\
	-d_{m} &
	\begin{array}
	[c]{cc}%
	\dotsc & -d_{1} \label{ABCD}
	\end{array}
	\end{array}
	\right]  ,\text{ }\mathbf{B}=\left[
	\begin{array}
	[c]{c}%
	\mathbf{0}\\
	1
	\end{array}
	\right]  ,\\
	\mathbf{C}  &  = \left[n_{m}, \ldots,n_{1}\right]. \nonumber
	\end{align}
\end{subequations}

In (\ref{ABCD}) $\mathbf{I}_{m-1}$ denotes the $m-1$ dimensional
identity matrix, and $\mathbf{0}$ denotes a column vector of zeros of length $m-1$.
The system matrix $\mathbf{A}$ is non-singular because $d_{m}\ne0$.

There is a system closely related to the estimation system that we refer to
as the dual system. Its input and output sequences are $(y_{j}^{\ast})_{j=1}^{\infty}$
and $(u_{j}^{\ast})_{j=1}^{\infty},$ and the dual plant, denoted
$P^{\ast},$ has the transfer function representation%

\begin{equation}
P^{\ast}(\lambda)=-\hat{\mathbf{n}}_{\mathrm{dual}} (\lambda)/ \hat{\mathbf{d}}_{\mathrm{dual}}(\lambda) 
\label{regconv}%
\end{equation}
where $\mathbf{n}_{\mathrm{dual}}=\left(  n_{m+1},\ldots,n_{1}\right)  $ and
$\mathbf{d}_{\mathrm{dual}}=\left(  d_{m+1},\ldots,d_2,1\right)$.
A minimal
state-space realization of the dual system is%

\begin{align}
\mathbf{x}_{k+1}^{\ast}  &  =\mathbf{A}^{\ast}\mathbf{x}_{k}^{\ast}+\mathbf{B}^{\ast}y_{k}%
^{\ast}\label{adjointstate1}\\
u_{k}^{\ast}  &  =\mathbf{C}^{\ast}\mathbf{x}_{k}^{\ast}+D^{\ast}y_{k}^{\ast}
\label{adjointstate}%
\end{align}%
\begin{align}
\mathbf{A}^{\ast}  &  =\left[
\begin{array}
[c]{cc}%
-d_{m-1}/d_{m} & \mathbf{I}_{m-1}\\%
\begin{array}
[c]{c}%
\vdots\\
-1/d_{m}%
\end{array}
&
\begin{array}
[c]{c}%
\\
\mathbf{0}
\end{array}
\end{array}
\right]  ,\label{ABCDstardef}\\
\text{ }\mathbf{B}^{\ast}  &  =\left[
\begin{array}
[c]{c}%
n_{m-1}\\
\vdots\\
n_{0}%
\end{array}
\right]  -\left[
\begin{array}
[c]{c}%
d_{m-1}\\
\vdots\\
1
\end{array}
\right]  \frac{n_{m}}{d_{m}},\label{ABCDstardef1}\\
\mathbf{C}^{\ast}  &  =\left[
\begin{array}
[c]{cccc}%
-1/d_{m} & 0 & \dotsc & 0
\end{array}
\right]  ,\text{ }D^{\ast}=\frac{-n_{m}}{d_{m}}.\label{ABCDstardef2}\\\nonumber
\end{align}
 
 \subsection{Polytopes}
 The primal and dual states defined in the previous Section will be interpreted in terms of the geometry of the polytopic uncertainty set $S_k$, so in this Section we introduce notation and briefly summarise the relevant theory of convex polytopes. For more information and background, including the definition of a polytope, see for example 
 \cite{Lay_convsets} or \cite{ziegler_polytopes}.
 Let $S$ denote a closed, convex polytope.
 
 The \textit{support function} of 
 $S$ in $\mathbb{R}^{m}$ is 
 \[h_{S}(\mathbf{f}%
 )=\max_{\mathbf{x}\in S}\left\langle \mathbf{f},\mathbf{x}\right\rangle
 ,\] where $\mathbf{f}\in\mathbb{R}^{m}.$
 
 When $\mathbf{f} \ne \mathbf{0}$ the set
 \[
 H_{S}(\mathbf{f}):=\left\{  \mathbf{x}\in\mathbb{R}^{m}:\left\langle
 \mathbf{f},\mathbf{x}\right\rangle =h_{S}(\mathbf{f})\right\}
 \]
 is the supporting hyperplane of $S$ with direction (outer normal vector)
 $\mathbf{f}.$  If $\mathbf{f}=\mathbf{0}$ then $H_{S}(\mathbf{f})=\mathbb{R}^{m}$.
 
 The intersection of $S$ with a supporting hyperplane is
 called a \textit{face} of $S,$ and a face of dimension $m-1$ is a \textit{facet} of $S$.  A face of dimension $m-2$ is called a \textit{ridge}, and the faces of dimensions $0$ and $1$ are termed \textit{vertices} and \textit{edges}, respectively.
  
 The \textit{normal cone} at a boundary point $\mathbf{x}$ is the set \[\left\{
 \mathbf{f}\in\mathbb{R}^{m}:\left\langle \mathbf{f},\mathbf{x}\right\rangle
 =h_{S}(\mathbf{f})\right\}  \]
 and is denoted $\mathcal{N}_S\left( \mathbf{x}\right)$.
 It is generated by the outward normals to the facets that form the polytope at $\mathbf{x}$, that is
 \[
 \mathcal{N}_{S}\left( \mathbf{x}\right)=\{\lambda_1\mathbf{f}_1+\dotsc+\lambda_n \mathbf{f}_n:\lambda_1,\ldots,\lambda_n\geq0\},
 \]
 where $\mathbf{f}_1,\ldots,\mathbf{f}_n$ are the directions of the facets containing $\mathbf{x}$. 
 Thus $\mathcal{N_S}\left(
 \mathbf{x}\right)$ contains the directions of all
 hyperplanes which touch $S$\ at $\mathbf{x.}$ 
 
 The sets of vertices and boundary points of $S$ are respectively denoted $\operatorname*{vert}(S)$ and $\partial S$. The interior of $S$ is denoted $\operatorname*{int}\left(S\right)$.
   If $\mathbf{x} \in \operatorname*{int}\left(S\right)$, then $\mathcal{N_S}\left(%
 \mathbf{x}\right):=\left\{\mathbf{0}\right\}$. By definition the directions of facets, and of the hyperplanes that contain facets, point outwards from the polytope. A direction is a non-zero vector. 
 
 The dual state $\mathbf{x}^*_k\left(\mathbf{y}^*,\mathbf{u}^*\right)$ can be interpreted as a direction vector $\mathbf{f}$ in the normal cone of the primal state $\mathbf{x}_k\left(\mathbf{y},\mathbf{u}\right) \in S_k$. We shall sometimes drop the subscript $k$, and indicate the next time instant with the superscript $^+$. For example, $\mathbf{f}$ will denote $\mathbf{x}^*_{k}$, and $\mathbf{f^+}$ will denote $\mathbf{x}^*_{k+1}$. Likewise, states $\mathbf{x}$ and facets $F$  belong to $S_k$, while $\mathbf{x}^+$ and $F^+$ belong to $S_{k+1}$.

\subsubsection{Geometric and combinatorial description of $S_{k}$} \label{Sectdescriptunc}

Three matrices, $\mathbf{V}$, $\mathbf{F}$ and $\mathcal{I}$ are used to describe the polytope $S_{k}$, and the algorithm updates all of them. Geometric data are contained in $\mathbf{V}$ and $\mathbf{F}$, while combinatorial information on the structure of $S_{k}$ is contained in $\mathcal{I}$. 

\begin{enumerate}
	\item   The $m$ by $n_v$ \textit{vertex matrix} $\mathbf{V}$ has columns containing the coordinates of the
	vertices of $S_{k}$.
	
	\item   The $m$ by $n_f$ \textit{facet matrix} $\mathbf{F}$ has columns containing the coordinates of the
	directions of the facets of $S_{k}$. Thus $\mathbf{F} = \left[\mathbf{f}_1,\dotsc,\mathbf{f}_{n_f}\right],$ where $\mathbf{f}_i$ is a column vector representing the direction of facet number $i.$ The magnitude of $\mathbf{f}_i$ is not important, that is, for any $\alpha >0$,  $\alpha \mathbf{f}_i$ can be used in place of $\mathbf{f}_i.$ 
	\item The \textit{vertex-facet incidence matrix} of $S_{k}$ is the matrix $\mathcal{I} \in \{0,1\}^{n_f \times n_v}$ which has entry 
	$\mathcal{I}(i,j)=1$ if $\mathbf{v}_j \in F_i$, and $\mathcal{I}(i,j)=0$ otherwise. So the facet $F_i$ is the convex hull of the vertices identified by the ones in the $i^{\rm th}$ row of $\mathcal{I}$.
\end{enumerate}

\section{Statement of Required Theorems}
 
Following Witsenhausen, \cite{Witsenhausen-poss-states-1968}, $S_{k+1}$ is given
recursively in terms of $S_{k}$ and the new observation $z_{k}$ by

\begin{equation} 
S_{k+1}=\left\{ \mathbf{x}^+ \; : \;
\begin{array}
[c]{l}%
\mathbf{x}\in S_{k}, \mathbf{x}^+=\mathbf{A}\mathbf{x}+\mathbf{B}u_{k},\\
y_{k}=\mathbf{C}\mathbf{x},\\
\left\vert u_{k}\right\vert \leq1,\text{ }\left\vert y_{k}-z_{k}\right\vert
\leq1.
\end{array} 
\right\}   \label{wits_recurs}%
\end{equation}

Special terminology is now introduced to describe states $\mathbf{x}$ and $\mathbf{x}^+$ related as in (\ref{wits_recurs}).

\begin{definition}
\label{defprecursor}The state $\mathbf{x} \in S_{k}$ is said to be a
\textit{precursor} of the state $\mathbf{x}^+$, $\mathbf{x}$
is \textit{propagated} to $\mathbf{x}^+,$ and $\mathbf{x}^+$ is a \textit{successor} to
$\mathbf{x},$ if there exists a scalar $u_k$ satisfying (\ref{primalstate1}), (\ref{primalstate2}), $\left\vert
u_k\right\vert \leq1$, and $\left\vert y_k-z_{k}%
\right\vert \leq1$.
\end{definition}

So $S_{k+1}$ is the set of all successors to all states in $S_{k}$, and any precursor of any state $\mathbf{x}^+\in S_{k+1}$ is in $S_{k}$.

We now define a relation between the inputs and outputs of the primal and dual systems.

\begin{definition}
	\label{def_align}The scalar pair $\left(  y_k,u_k\right)  $ is said to be \textit{aligned}
	with $\left(  y_k^{\ast},u_k^{\ast}\right)  $ if
	\begin{equation}%
		\begin{array}
			[c]{c}%
			u_k^{\ast}>0\implies u_k=1\\
			u_k^{\ast}<0\implies u_k=-1\\
			\left\vert u_k\right\vert <1\implies u_k^{\ast}=0
		\end{array}
		\label{align_defv}%
	\end{equation}
	and%
	\begin{equation}%
		\begin{array}
			[c]{c}%
			y_k^{\ast}>0\implies y_k=1+z_{k}\\
			y_k^{\ast}<0\implies y_k=-1+z_{k}\\
			\left\vert y_k-z_{k}\right\vert <1\implies y_k^{\ast}=0.
		\end{array}
		\label{align_defy}%
	\end{equation}
	
\end{definition}

Associated with any $\mathbf{f}\in \mathbb{R}^{m}$ define in the
$y^{\ast}_{k}u^{\ast}_{k}$-plane the dual line $L^{\ast}\left(\mathbf{f}\right)$:

\begin{definition}
	$L^{\ast}\left(\mathbf{f}\right)=\{(y^*_k,u^*_k):n_{m}y^{\ast}_{k}+d_{m}u^{\ast}_{k}=-\left(\mathbf{f}\right)_1\}.$
\end{definition}

The scalars $y^*_k$ and $u^*_k$, the input and output of the dual plant at time $k$, are constrained to lie on the line $L^{\ast}\left(  \mathbf{f}\right)$ by (\ref{adjointstate}) and (\ref{ABCDstardef2}). 

The square $Q$ in the $y_ku_k$-plane contains points $(y_k,u_k)$ satisfying $\left\vert u_{k}\right\vert \leq1$ and $\left\vert y_{k}-z_{k}\right\vert
\leq1$.
 
One more definition is required in order to state the Theorems quoted from \cite{Hill-luo-Schw-2016}.
\begin{definition}
	\label{defM}Given $\mathbf{x}\in S_{k}$, $\mathbf{f}\in \mathcal{N}_{S_{k}}(\mathbf{x})$  and $z_{k},$ the set $M\left(
	\mathbf{x},\mathbf{f},z_{k}\right)$ is the set of scalar pairs $\left(  u_k,y_k^{\ast}\right)$ which satisfy
	\begin{enumerate}
		\item $	(y_k,u_k) \in Q$, where  $y_k = \mathbf{C}\mathbf{x}$ and
		\item  $\left(y_k,u_k\right)$   is aligned with $\left(
	y_k^{\ast},u_k^{\ast}\right)$,  where
	$(y_k^{\ast},u_k^{\ast}) \in L^{\ast}\left(  \mathbf{f}\right).$
   \end{enumerate}
\end{definition}

\begin{figure}[ptb]
	\begin{center} 
		\includegraphics[scale=0.5]{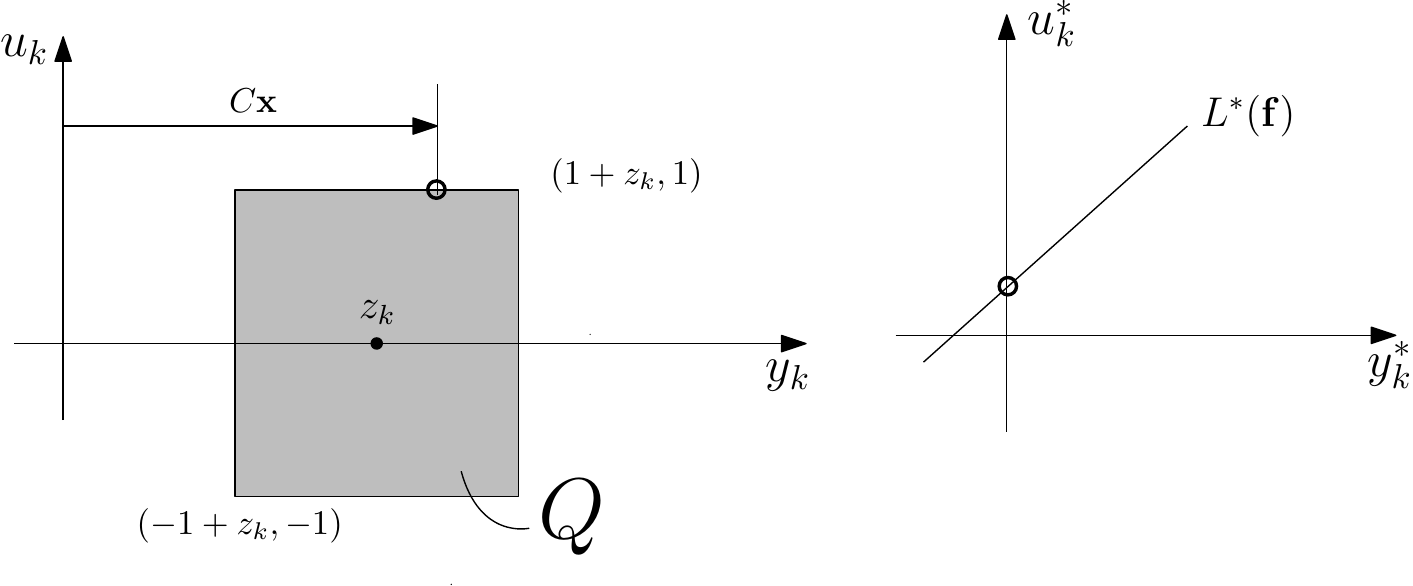} \caption{The circled points, at $\left(C\mathbf{x},1\right)$ in the $y_ku_k$ plane, and at the intersection of $L^{\ast}(\mathbf{f})$ and the $u^*_k$ axis, are aligned, so $M\left(
			\mathbf{x},\mathbf{f},z_{k}\right)=\{(1,0)\}$ because at the aligned points $u_k=1$ and $y^*_k=0$.} \label{Q_L_intersect}%
	\end{center}
\end{figure}

 Finding $M$ is computationally very simple. For example, in Fig. \ref{Q_L_intersect} alignment occurs solely between the two circled points, so $M$ contains the singleton pair $(u_k,y^*_k)=(1,0)$.

\begin{figure*}[ptb]
	\centering
	\includegraphics[scale=0.9]{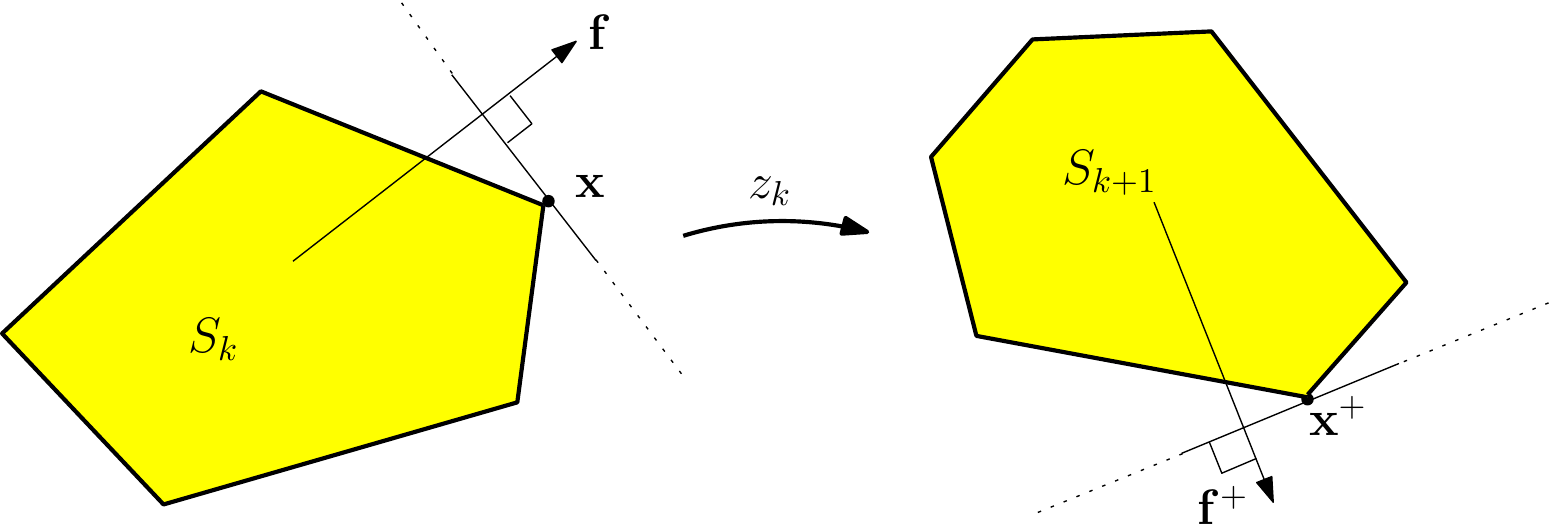} \caption{The vectors in Theorem \ref{main1}. The state
		$\mathbf{x}\in\partial S_{k}$ and direction $\mathbf{f}\in
		\mathcal{N}_{S_{k}}(\mathbf{x})$ are propagated to $\mathbf{x^+}%
		\in\partial S_{k+1}$ and $\mathbf{f^+}\in \mathcal{N}_{S_{k+1}}(\mathbf{x^+}%
		)$ by the measurement $z_{k}.$}%
	\label{basic_conv}%
\end{figure*}

We now present the Theorems from which the algorithm is derived. See Fig. \ref{basic_conv} for a geometric depiction of the vectors in Theorem \ref{main1} for the important special case where $\mathbf{x}$ and $\mathbf{x}^+$ are both boundary points. Proofs are in \cite{Hill-luo-Schw-2016}.

\begin{theorem}
	\label{main1}Suppose $\mathbf{x}\in S_{k}$ and
	$\mathbf{f}\in \mathcal{N}_{S_{k}}(\mathbf{x})$. Then $\mathbf{x^+}=\mathbf{A}\mathbf{x}+\mathbf{B}u_{k} \in S_{k+1}$ and 
	$\mathbf{f^+}
	=\mathbf{A}^{\ast}\mathbf{f}+\mathbf{B}^{\ast}y_{k}^{\ast
	}\in \mathcal{N}_{S_{k+1}}(\mathbf{x^+})$ if and only if $\left(
	u_{k},y_{k}^{\ast}\right)  \in M\left( \mathbf{x}, \mathbf{f},z_{k}\right)$.
\end{theorem}

 In a typical application  a point on the boundary of $S_{k}$, and the direction of any one of its supporting hyperplanes, are propagated to a point on the boundary of $S_{k+1}$ along with the direction of one of its supporting hyperplanes.  

The companion result Theorem \ref{main2}, given below, shows that any boundary point $\mathbf{x^+} \in S_{k+1}$, and any direction in the normal cone of $\mathbf{x^+}$, are attainable from \textit{any} precursor $\mathbf{x}$ of $\mathbf{x^+}$  and \textit{some} direction in the normal cone of $\mathbf{x}$. 

\begin{theorem}
	\label{main2}	
	Select any $\mathbf{x^+}\in S_{k+1}$, any $\mathbf{f^+} \in \mathcal{N}_{S_{k+1}}\left( \mathbf{x^+}\right)$ and any precursor $\mathbf{x}$ of $\mathbf{x^+}$. There exists $\mathbf{f} \in \mathcal{N}_{S_{k}}\left( \mathbf{x}\right)$ and $\left(
	u_{k},y_{k}^{\ast}\right)  \in M\left( \mathbf{x},  \mathbf{f},z_{k}\right)$
	for which $\mathbf{x^+} =\mathbf{A}\mathbf{x}+\mathbf{B}u_{k}$ and $\mathbf{f^+}=\mathbf{A}^{\ast}\mathbf{f}+\mathbf{B}^{\ast}y_{k}^{\ast}$.
\end{theorem}

In the rest of the paper we show how to use the propagation of states and their normal cones to update the uncertainty set for plants with a lag.  An implementation of our algorithm in Matlab code, \texttt{uncertaintysetwithlag.m}, is available on the link below\footnote{Go to \textrm{https://mathworks.com/matlabcentral/fileexchange}}. Also given is code \texttt{tcomplag.m} which allows comparisons of computation time for our algorithm, denoted F-V (facets-vertices), Fourier-Motzkin  and mplp (multi-parametric linear programming). We used the software package \cite{mpt} for Fourier-Motzkin projection and mplp. See \cite{Shamma_Kuang-1999} for an explanation of the use of Fourier-Motzkin projection to compute updated uncertainty sets. 

\section{The incidence matrix for $S_k$}

We turn now to the construction of  $\mathcal{I}^+$, the incidence matrix for $S_{k+1}$, shown in Table \ref{incidmatrix}, where the six incidence matrices in the body of the Table are defined in Section \ref{sect_table}. Initially $\mathcal{I}^+$ is put equal to the $\left(3n_f + n_R \right)$ by $2n_v$ zero matrix, that is $\mathcal{I}^+$ is allocated sufficient space to accommodate all conceivable facets and vertices of $S_{k+1}.$ So there are $3n_f + n_R$ rows in the Table, where $n_f$ and $n_R$ are respectively the number of facets and ridges in $S_k$. The rows are divided into fours block, each of the first three blocks comprising $n_f$ rows, and the last block $ n_R$ rows. The first two blocks point to facets of $S_{k+1}$ whose directions have the property  $\left(\mathbf{f}\right)_m \ne 0$, while the last two blocks point to facets whose directions $\mathbf{f}$ have the property  $\left(\mathbf{f}\right)_m = 0$.

Turning now to the columns of Table \ref{incidmatrix}, the first block of $n_v$ columns points to states $\mathbf{A}\mathbf{v}_j+\mathbf{B}$,  $\mathbf{v}_l \in \operatorname{vert}\bigl(S_k\bigr)$, which, depending on $z_k$, might be vertices of $S_{k+1}$. The second block of $n_v$ columns points to states $\mathbf{A}\mathbf{v}_j-\mathbf{B}$, also possible vertices of $S_{k+1}$.

After $\mathcal{I}^+$ has been constructed in accordance with Table \ref{incidmatrix} it must have some rows and columns composed entirely of zeros. The final step is for these rows and columns of zeros to be deleted.

\begin{table}
	\begin{center}
		\begin{tabular}
			{|l|c|c|} \hline
			\diaghead{\theadfont ColumnmnHead }%
			{Facet \\ directions}{Vertices} & 
			\thead{  $g^{\mathrm{T}}\left(\mathbf{V}\right)$  } & \thead{$g^{\mathrm{B}}\left(\mathbf{V}\right)$ } \\    \hline
			\rule{0pt}{2.5ex}$A^*\mathbf{f}_i$(with $u_k = 1$)
			& $\mathcal{I}^{\mathrm{T}}$ & $\mathbf{0}$   \\   \hline 
			\rule{0pt}{2.5ex} $A^*\mathbf{f}_i$(with $u_k = -1$) 
			& $\mathbf{0}$ & $\mathcal{I}^{\mathrm{B}}$ \\ 
			\hline
			\rule{0pt}{3ex}  $A^*\mathbf{f}_i$(where $\left(\mathbf{f}\right)_1=0$)
			& $\mathcal{I}O^{\mathrm{T}}$ & $\mathcal{I}O^{\mathrm{B}}$  \\    \hline
			\rule{0pt}{3ex} $A^*\mathbf{f}^R_l$(Facets of $S_{k+1}$ from ridges)
			& $\mathcal{I}R^{\mathrm{T}}$ & $\mathcal{I}R^{\mathrm{B}}$  \\    \hline
		\end{tabular}
	\end{center}
	\caption{The incidence matrix for $S_{k+1}$ prior to the removal of rows and columns composed entirely of zeros}%
	\label{incidmatrix}%
\end{table}

\section{Filling in Table \ref{incidmatrix}} \label{sect_table}

We introduce notation to describe two ways a state can be propagated from  $S_{k}$ to $S_{k+1}$. 
The mappings $g^\mathrm{T}$ and $g^\mathrm{B}$, where the superscripts indicate the top and bottom sides of $Q$, map a state  on the boundary of $S_{k}$ to the boundary of $S_{k+1}$.

On the top (bottom) side of $Q$, $u_k=1$ $(u_k=-1)$, and (\ref{primalstate1})
is used to propagate $\mathbf{x}$. 

\begin{definition}
	For any $\mathbf{x} \in S_{k}$ put
	\begin{equation}
	\begin{array}
	[c]{c}%
	g^\mathrm{T}\left(\mathbf{x}\right) = \mathbf{A}\mathbf{x} +\mathbf{B}\\
	g^\mathrm{B}\left(\mathbf{x}\right) = \mathbf{A}\mathbf{x} -\mathbf{B}.\\
	\label{defsidesQprop}%
	\end{array}
	\end{equation}
\end{definition}

\begin{definition} \label{defpropup}
	A point $\mathbf{x} \in S_{k}$ is said to \textit{propagate up} (resp. \textit{down}) if  $g^\mathrm{T}\left(\mathbf{x}\right)$ $\bigl( (g^\mathrm{B}\left(\mathbf{x}\right)\bigr)$ lies on the boundary of $S_{k+1}$.
\end{definition} 

These propagations occur when the maximum or minimum allowable value of the input ($\pm$1) is applied to the plant, and furthermore the resulting state of the plant lies on the boundary of $S_{k+1}$. Theorem \ref{main1} tells us, given $\mathbf{x}$ and $\mathcal{N}_{S_{k}}\left( \mathbf{x}\right)$, if $\mathbf{x}$ propagates up or down. For example, in Fig. \ref{Q_L_intersect}, if $\mathbf{x} \in \partial S_k$ then  $\mathbf{x}$ propagates up to $g^\mathrm{T}\left(\mathbf{x}\right)$ because there exists $\mathbf{0} \ne \mathbf{f} \in \mathcal{N}_{S_{k}}\left( \mathbf{x}\right)$ which implies $(u_k,y^*_k) = (1,0) \in  M\left( \mathbf{x}, \mathbf{f},z_{k}\right)$ and so, by Theorem \ref{main1}, $\mathbf{x}^+=\mathbf{A}\mathbf{x}+\mathbf{B} \in S_{k+1}$ and $\mathbf{f^+}
=\mathbf{A}^{\ast}\mathbf{f}\in \mathcal{N}_{S_{k+1}}(\mathbf{x^+})$ and is non-zero, further implying that $\mathbf{x}^+ \in \partial S_{k+1}$.

Definition \ref{defpropup} can be extended to facets as follows.
\begin{definition} \label{deffacetspropup}Suppose $F$ is a facet $S_k$, and $F^+$ is a facet $S_{k+1}$. Then $F$ is said to propagate up (down) to $F^+$ if $F^+ = g^T(F)$, $ \bigl(F^+ = g^B(F)\bigr)$.
\end{definition}

Here $g^{\mathrm{T}}\left(F\right)$, for example, denotes the set $\left\{g^{\mathrm{T}}\left(\mathbf{x}\right):\mathbf{x} \in F\right\}$.

If $F$ propagates either up or down to $F^+$ then $F^+$ is an affinely isomorphic image of $F$ because $\mathbf{A}$ is invertible. 

Let $F_i$ be any facet of $S_k$, and denote its direction by $\mathbf{f}_i$. There are three possibilities:
\begin{enumerate}
	\item The line $L^*\bigl(\mathbf{f}_i\bigr)$ intersects the positive $u^*_k$ axis, that is $d_m\bigl(\mathbf{f}_i\bigr)_1 < 0$. The $n_f \times 1$ logical vector $\mathcal{I}F^{\mathrm{T}}$ points to facets with such directions. Explicitly,  $\mathcal{I}F^{\mathrm{T}}(i) =1$ if $\mathbf{f}_i$ satisfies $d_m\bigl(\mathbf{f}_i\bigr)_1 < 0$, and $\mathcal{I}F^{\mathrm{T}}(i) =0$ otherwise; 
	\item  $L^*\bigl(\mathbf{f}_i\bigr)$ intersects the negative $u^*_k$ axis, that is $d_m\bigl(\mathbf{f}_i\bigr)_1 >0$.  Define $\mathcal{I}F^{\mathrm{B}}(i) =1$ if $\mathbf{f}_i$ satisfies $d_m\bigl(\mathbf{f}_i\bigr)_1 > 0$, and $\mathcal{I}F^{\mathrm{B}}(i) =0$ otherwise; and
	\item $L^*\bigl(\mathbf{f}_i\bigr)$ passes through the origin. Define $\mathcal{I}F^{\mathrm{O}}(i) =1$ if $\mathbf{f}_i$ satisfies $\bigl(\mathbf{f}_i\bigr)_1 = 0$, and $\mathcal{I}F^{\mathrm{O}}(i) =0$ otherwise.
\end{enumerate}
\subsubsection{Case 1}For every $\mathbf{x}_j$ in $F_i$ we have the situation depicted in Fig. \ref{Q_L_intersect}. By Theorem \ref{main1}, $\mathbf{x}_j$ is propagated up, and $F_i$ is propagated up to the facet  $F_i^+:=g^T\bigl(F_i\bigr)$ of  $S_{k+1}$, which has direction $\mathbf{A}^*\mathbf{f}_i$.
The logical matrix $\mathcal{I}^{\mathrm{T}} \in \{0,1\}^{n_f \times n_v}$ is defined to have entry 
$\mathcal{I}^{\mathrm{T}}(i,j)=1$ if
 $\mathcal{I}(i,j)=1$ and
	$\mathcal{I}F^{\mathrm{T}}(i)=1$,
and  $\mathcal{I}^{\mathrm{T}} (i,j)=0$ otherwise.
A one in the  $i^{\rm th}$ row and  $j^{\rm th}$ column of $\mathcal{I}^{\mathrm{T}}$ implies that there is a facet, $F_i^+$, of $S_{k+1}$ with direction $\mathbf{A}^*\mathbf{f}_i$, and that $g^\mathrm{T}\left(\mathbf{v}_j\right) \in F_i^+$.

Case 2) is similar to Case 1), except that $F_i$ is propagated down. The matrix $\mathcal{I}^{\mathrm{B}} \in \{0,1\}^{n_f \times n_v}$ is defined to have entry 
$\mathcal{I}^{\mathrm{B}}(i,j)=1$ if
$\mathcal{I}(i,j)=1$ and
$\mathcal{I}F^{\mathrm{B}}(i)=1$,
and  $\mathcal{I}^{\mathrm{B}} (i,j)=0$ otherwise.

We require for Case 3), covered below, as well as for the case of propagation of facets from ridges, to be treated in the next section, two logical matrices which point to those vertices of $S_k$ that are propagated to vertices of $S_{k+1}$. The $ n_v \times 1$ vector $\mathcal{I}V^{\mathrm{PT}}$ is  defined by $\mathcal{I}V^{\mathrm{PT}}(j)=1$ if the $j^{\rm th}$ column of $\mathcal{I}^{\mathrm{T}}$ is not composed entirely of zeros, and $\mathcal{I}V^{\mathrm{PT}}(j)=0$ otherwise. Thus $\mathcal{I}V^{\mathrm{PT}}$ points to those vertices which propagate up to a vertex of $S_{k+1}$. Likewise $\mathcal{I}V^{\mathrm{PB}}$ is defined by $\mathcal{I}V^{\mathrm{PB}}(j)=1$ if the $j^{\rm th}$ column of $\mathcal{I}^{\mathrm{B}}$ is not composed entirely of zeros, and $\mathcal{I}V^{\mathrm{PB}}(j)=0$ otherwise.

Every point of $Q$ is aligned with the origin in the $y^*_ku^*_k$ plane so, in Case 3), by Theorem \ref{main1} every $\mathbf{x}_j$ in $F_i$ is propagated to $\mathbf{x}_j^+ =\mathbf{A}\mathbf{x}_j+\mathbf{B}u_k \in \partial S_{k+1}$ for all $\left\vert u_{k}\right\vert \leq1$. 
 We show the method for constructing all vertices of the facet $F_i^+$, and provide a sketch of the proof that it works. Put 
$F_i^+=\{\mathbf{x}_j^+: \mathbf{x}_j^+ =\mathbf{A}\mathbf{x}_j+\mathbf{B}u_k \text{ where } \mathbf{x}_j \in F_i, u_k \in [-1,1]\}.$
A convexity argument shows
\[ F_i^+=\{\mathbf{x}_j^+: \mathbf{x}_j^+ =g^\mathrm{T}\left(\mathbf{v}_j\right)\text{ or } \mathbf{x}_j^+ =g^\mathrm{B}\left(\mathbf{v}_j\right), \mathbf{v}_j \in \operatorname{vert}(F_i) \}. \]
 We seek the vertices of the facet $F_i^+$. Now every precursor $\mathbf{v}_j$ of any $\mathbf{v}_j^+$ in the vertex set of $F_i^+$ is a vertex of $F_i$. However, it is not necessarily the case that every successor of every vertex of $F_i$ is a vertex of $F_i^+$. We must remove from $F_i^+$ those $\mathbf{x}_j^+$ which are not vertices of $F_i^+$, and this can be done by removing those $\mathbf{x}_j^+$ which are not vertices of $S_k$.
 
Define the incidence matrix $\mathcal{I}O^{\mathrm{T}}
(i,j)=1$ if $\mathcal{I}F^{\mathrm{O}}(i) =1$, $\mathcal{I}V^{\mathrm{PT}}(j)=1$  and $\mathcal{I}(i,j)=1$, and $\mathcal{I}O^{\mathrm{T}}
(i,j)=0$ otherwise. In words, $\mathcal{I}O^{\mathrm{T}}(i,j)=1$ if $\mathbf{v}_j$ propagates up to a vertex of $S_k$, and $\mathbf{v}_j$ 
 is in $F_i$, with $(\mathbf{f}_i)_1=0$. 
 
 The matrix $\mathcal{I}O^{\mathrm{B}}$, pointing to those vertices of facets with direction vector whose first component is zero that propagate down, is defined similarly. Thus $\mathcal{I}O^{\mathrm{B}}
 (i,j)=1$ if $\mathcal{I}F^{\mathrm{O}}(i) =1$, $\mathcal{I}V^{\mathrm{PB}}(j)=1$  and $\mathcal{I}(i,j)=1$, and $\mathcal{I}O^{\mathrm{B}}
 (i,j)=0$ otherwise.
 
 This completes the construction of all of Table \ref{incidmatrix} except the last block.

\subsection{Propagation from a ridge in $S_{k}$ to a facet in $S_{k+1}$} \label{sectridgetofacet}
In the previous Section Case 3) was more intricate than Cases 1) and 2) because it treated facets with directions whose dual lines passed through the origin, opening up an infinite number of alignment possibilities. There is one, and only one, other situation in which a boundary point $\mathbf{x}$ of $S_k$ can have a direction $\mathbf{f}$ in its normal cone which satisfies $\bigl(\mathbf{f}\bigr)_1 = 0$. If such $\mathbf{x}$ does not lie on a facet covered by case 3) above,  $\mathbf{x}$ must belong to a special type of ridge, to be defined next.

\begin{figure}[!h]
	\centering
	\includegraphics[scale=0.7]{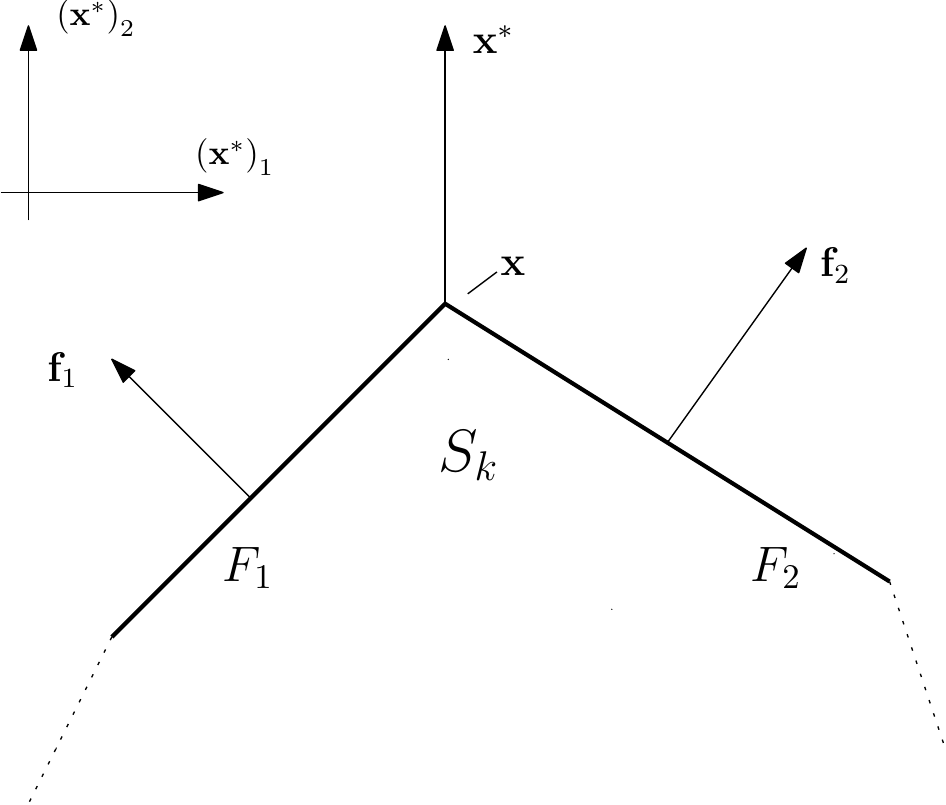} \caption{Two facets intersect in a ridge containing $\mathbf{x}$. The direction $\mathbf{x}^{\ast} \in \mathcal{N}_{S_{k}}(\mathbf{x})$ satisfies $\left( \mathbf{x}^{\ast}\right) _1 = 0$.}
	\label{ridge}%
\end{figure}

See Fig. \ref{ridge} where $S_{k}$ possesses two intersecting facets, with directions $\mathbf{f}_1$ and $\mathbf{f}_{2}$, whose first components have opposite signs. These facets have the state $\mathbf{x}$ in their intersection, so necessarily there exists a direction in the normal cone of $\mathbf{x}$ whose first component is zero. We define $R$ to be the set of all such ridges of $S_{k}$, and denote its cardinality by $n_R$. Clearly $R$ is non-empty if $S_{k}$ is non-empty and full-dimensional.

\begin{definition} 
	\label{defridges}
	Define $R$ to be the set of ridges, $R_l$, $l=1,\ldots,n_R$, of $S_k$ with the property that $R_l$ is the intersection of two facets $F_1$  and  $F_2$  say,  for which $\left( \mathbf{f}_1 \right)_1 <0$ and $\left(\mathbf{f}_2\right)_1 > 0$,  where  $\mathbf{f}_1$   and $\mathbf{f}_2$   are the directions of $F_1$  and  $F_2$.
\end{definition}


Associated with any ridge in $R$ there is by definition a direction in the cone generated by $\mathbf{f}_1$ and $\mathbf{f}_2$ having the property that its first component is zero. This direction is unique.

\begin{definition} \label{defridgedir} For any $R_l \in R$ define $\mathbf{f}_l^R$ to be the direction in the convex hull of $\mathbf{f}_1$ and $\mathbf{f}_2$ for which $\left(\mathbf{f}_l^R\right)_1 = 0$.
\end{definition}

Clearly $\mathbf{f}_l^R \in \mathcal{N}_{S_{k}}\left( \mathbf{x}\right)$ for all $\mathbf{x} \in R_l$.
The proof of the following Proposition is omitted.
\begin{proposition}
	For all $R_l \in R$ there is a facet of $S_k$ with direction $\mathbf{A}^*\mathbf{f}_l^R$.
\end{proposition}

For any $S_{k}$ and any $R_l \in R$, the direction $\mathbf{f}_l^R$ can be found using linear interpolation. More precisely $\mathbf{f}_l^R=t\mathbf{f}_1+(1-t)\mathbf{f}_2$ where $t=\left(\mathbf{f}_1\right)_1/[ \left(\mathbf{f}_2\right)_1 - \left(\mathbf{f}_1\right)_1 ]$, and $\mathbf{f}_1$ and $\mathbf{f}_2$ are the directions of the facets which intersect to form $R_l$. 

The incidence matrices in the bottom row of Table \ref{incidmatrix} can now be defined.

The vertex set of $R_l$ is the intersection of the vertex sets of the facets which intersect to form $R_l$. 
The matrix $\mathcal{I}RV \in \{0,1\}^{{n_{R}} \times n_v}$ has entry $\mathcal{I}RV(l,j)=1$ if the ridge $R_l$ contains $\mathbf{v}_j$, and $\mathcal{I}RV(l,j)=0$ otherwise. So the $l^{\rm th}$ row of $\mathcal{I}RV$ points to the vertex set of $R_l$. 

The matrix $\mathcal{I}R^{\mathrm{T}} \in \{0,1\}^{{n_{R}} \times n_v}$ has entry $\mathcal{I}R^{\mathrm{T}}(l,j)=1$ if $\mathcal{I}V^{\mathrm{PT}}(j)=1$ and $\mathcal{I}RV(l,j)=1$, and $\mathcal{I}R^{\mathrm{T}}(l,j)=0$ otherwise. So the $l^{\rm th}$ row of $\mathcal{I}R^{\mathrm{T}}$ points to vertices $\mathbf{v}_j$ of $R_l$ which propagate up to vertices of the ridge $R_l$ in $S_k$. 

The matrix $\mathcal{I}R^{\mathrm{B}} \in \{0,1\}^{{n_{R}} \times n_v}$ has entry $\mathcal{I}R^{\mathrm{B}}(l,j)=1$ if $\mathcal{I}V^{\mathrm{PB}}(j)=1$ and $\mathcal{I}RV(l,j)=1$, and $\mathcal{I}R^{\mathrm{T}}(l,j)=0$ otherwise. So the $l^{\rm th}$ row of $\mathcal{I}R^{\mathrm{B}}$ points to vertices $\mathbf{v}_j$ of $R_l$ which propagate down to vertices of the ridge $R_l$ in $S_k$. 

This completes the definition of all of the incidence matrices occurring in Table \ref{incidmatrix}. We claim that $\mathcal{I}^+$ constructed according to Table \ref{incidmatrix} is, after removal of rows and columns containing all zeros, the incidence matrix of $S_{k+1}$. We have sketched the proof, which ultimately relies on Theorems \ref{main1} and \ref{main2}.

\section{Examples}
Theorems \ref{main1} and \ref{main2} are most naturally interpreted as yielding $S_{k+1}$, consistent with measurements up to $z_k$. The input is $S_{k}$ consistent with measurements up to $z_{k-1}$. In common terminology we ``update then propagate". It is more common to determine $S_{k+1}$ consistent with measurements up to $z_{k+1}$, with input $S_{k}$ consistent with measurements up to $z_{k}$, essentially ``propagate then update". The accompanying code implements this latter case. For plants with a lag changing the order of the operations update and propagate is very easy.

\begin{figure}
	\centering
	\includegraphics[scale=0.6]{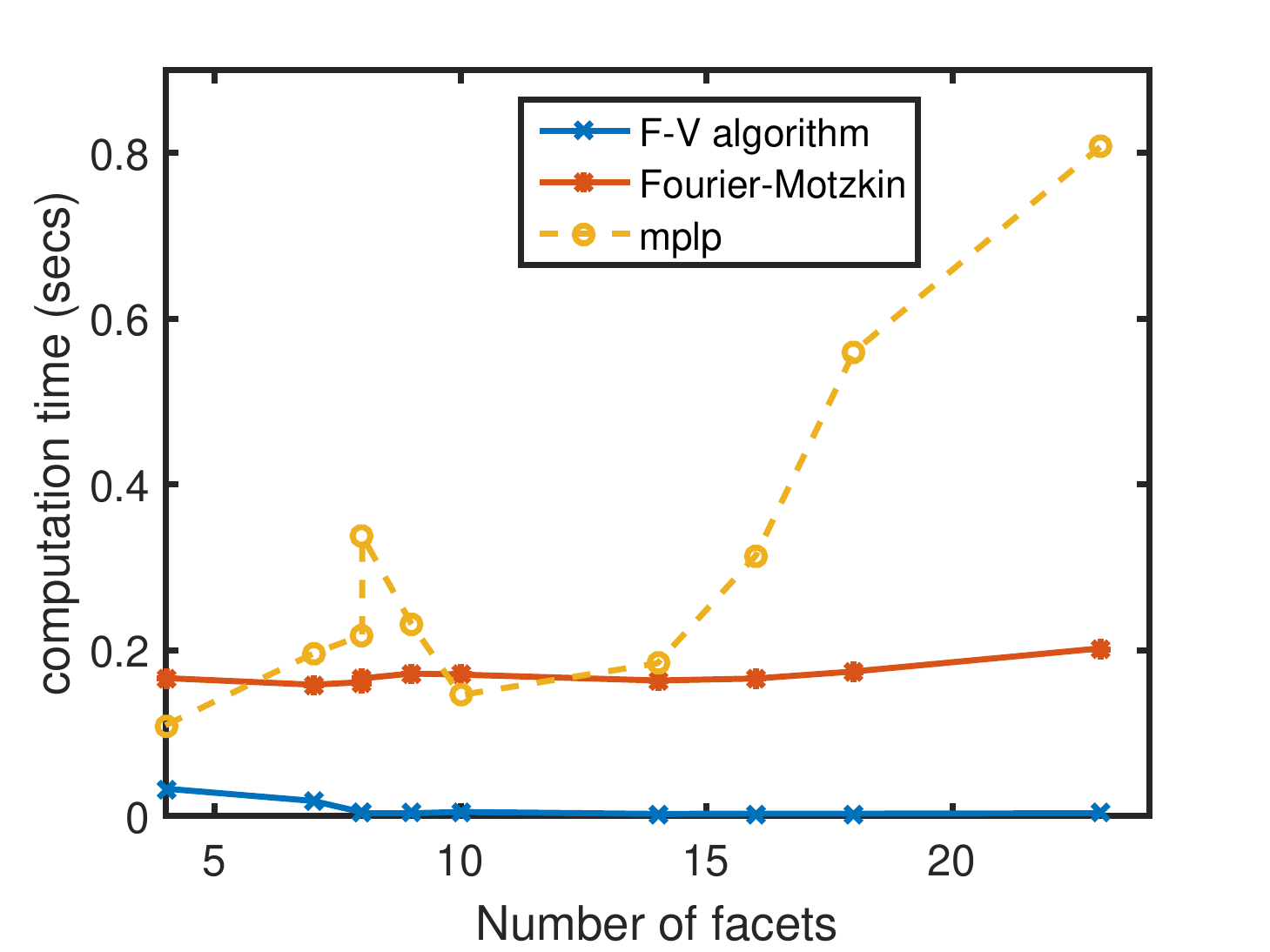} \caption{The time taken to update $S_k$ versus the number of facets of $S_k$ for a third order plant with a lag.}%
	\label{order3_kis10}%
\end{figure}

\begin{figure}
	\centering
	\includegraphics[scale=0.6]{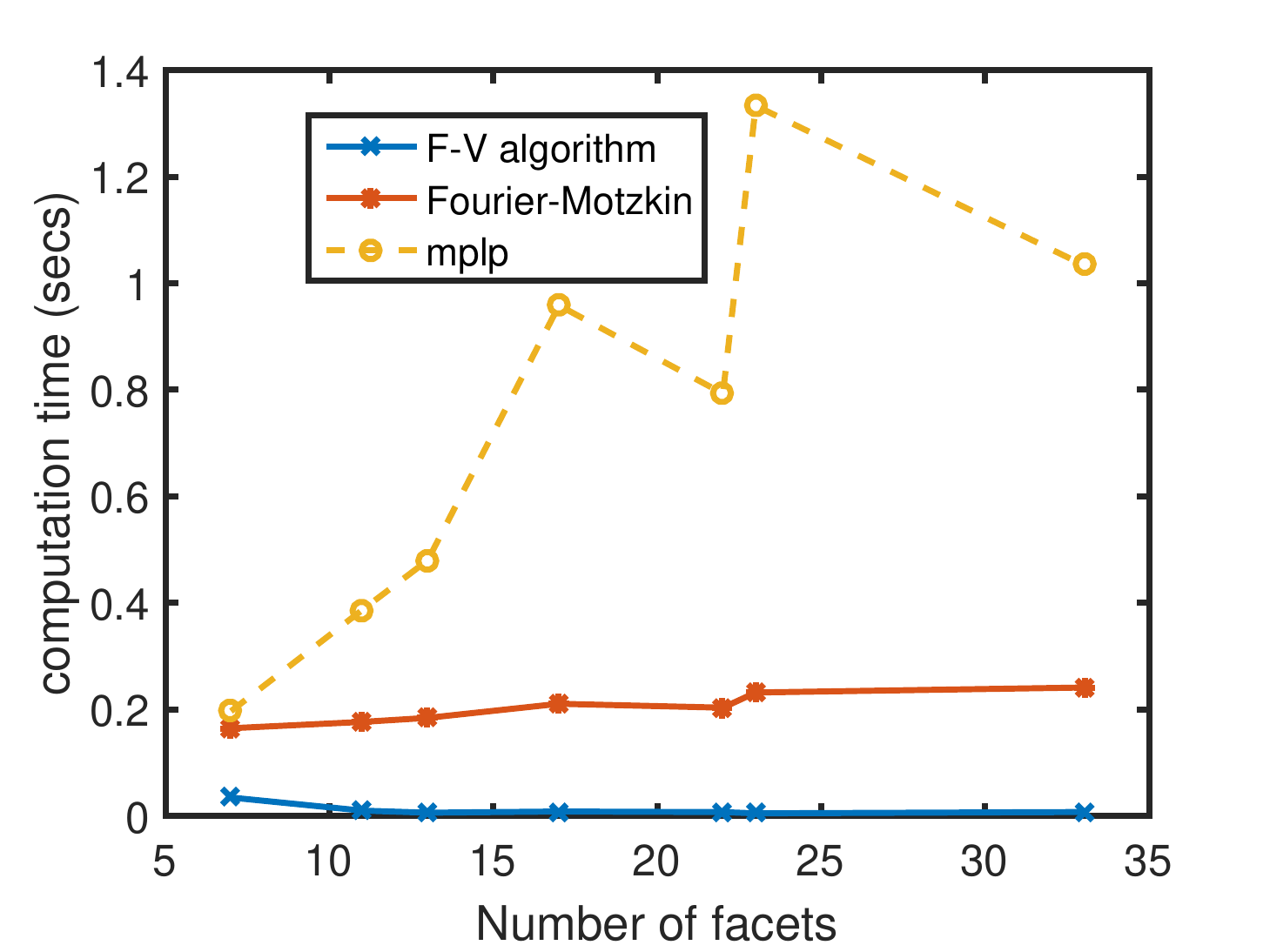} \caption{The time taken to update $S_k$ versus the number of facets of $S_k$ for a fifth order plant with a lag.}%
	\label{order5_kis7}%
\end{figure}

We give representative plots comparing our algorithm, denoted F-V, projection using Fourier-Motzkin elimination, and projection using multi-parametric linear programming. They were generated using the Matlab code accompanying the paper. The reader can generate plots for her/himself with differing plant sizes and lengths of measurement sequences.  Figure \ref{order3_kis10} compares the computation times of the three methods for a third order randomly selected stable plant, with randomly selected measurements. Figure \ref{order5_kis7} compares the computation times for a fifth order randomly selected stable plant, again with randomly chosen measurements. It can be seen that the algorithm presented in this paper is significantly faster. For these examples the F-V algorithm was between twenty and forty times faster, on average, than the better of the two other methods.

\section{CONCLUSIONS}
For the case of plants with a lag we have described a new method of uncertainty set propagation that appears to be more efficient than those in current use. There are in addition other advantages. It provides more information than most current methods, because the vertices and facets are both updated at every time step. Finally, because the computationally intensive steps are done using Boolean matrices, for moderately sized problems all calculations can be performed exactly over the rationals.

\bibliographystyle{siam}
\bibliography{acompat,Ref}

\end{document}